\documentclass[conference]{IEEEtran}
\usepackage{mystyle}

\usepackage[nolist]{acronym}
\begin{acronym}
	\acro{TSO}{transmission system operator}
	\acro{AC}{alternating current}
	\acro{HVDC}{high voltage direct current}
	\acro{API}{application programming interface}
	\acro{SCOPF}{security constrained optimal power flow}
	\acro{OPF}{optimal power flow}
    \acro{ED}{economic dispatch}
    \acro{IMML}{inverse matrix multiplication lemma}
	\acro{LP}{linear program}
    \acro{OPS}{optimal probabilistic security}
    \acro{SCOPF}{security-constrained optimal power flow}
    \acro{P-SCOPF}{preventive SCOPF}
    \acro{C-SCOPF}{corrective SCOPF}
    \acro{RMAC}{reliability management approach and criterion}
    \acro{IMML}{inverse matrix modification lemma}
    \acro{SM}{Sherman-Morrison}
    \acro{SMW}{Sherman-Morrison-Woodbury}
    \acro{DRC}{deterministic reliability criterion}
    \acro{MINLP}{mixed-integer nonlinear program}
    \acro{TSO}{transmission system operator}
    \acro{PTDF}{power transfer distribution factors}
    \acro{LODF}{line outage distribution factors}
    \acro{CENS}{cost of energy not served}
    \acro{ENS}{energy not served}
    \acro{LOL}{loss of load}
    \acro{FOR}{forced outage rate}
    \acro{EOC}{expected operational cost}
    \acro{FCR}{frequency containment reserves}
    \acro{aFRR}{automatic frequency restoration reserves}
    \acro{mFRR}{manual frequency restoration reserves}
    \acro{FFR}{fast frequency reserves}
    \acro{VOLL}{value of lost load}
    \acro{VRE}{variable renewable energy}
    \acro{RTS}{reliability test system}
\end{acronym}

\begin{document}
\title{A Fast and Scalable Iterative Solution of a Socio-Economic Security-Constrained Optimal Power Flow with Two-Stage Post-Contingency Control
\thanks{The research leading to results in this paper has received funding from the Research Council of Norway through the project ``Resilient and Probabilistic reliability management of the transmission grid'' (RaPid) (Grant No. 294754), The Norwegian Water Resources and Energy Directorate, and Statnett.}
}

\author{\IEEEauthorblockN{Matias Vistnes\IEEEauthorrefmark{1},
Vijay Venu Vadlamudi\IEEEauthorrefmark{1},
and Oddbjørn Gjerde\IEEEauthorrefmark{2}}
\IEEEauthorblockA{\IEEEauthorrefmark{1} Department of Electric Energy, 
Norwegian University of Science and Technology
}
\IEEEauthorblockA{\IEEEauthorrefmark{2} Department of Energy Systems, SINTEF Energy Research
}
Trondheim, Norway
}

\maketitle

\begin{abstract}
A fast and scalable iterative methodology for solving the security-constrained optimal power flow (SCOPF) problem is proposed using problem decomposition and the inverse matrix modification lemma. The SCOPF formulation tackles system security operational planning by using short- and long-term post-contingency limits, probability of branch outages, and preventive and corrective actions, a probabilistic corrective-SCOPF problem formulation. Using two post-contingency states and contingency probabilities, the SCOPF could provide good system security at a lower cost than only preventive actions as the typical `N-1'-formulation does. Additional security is ensured using a post-contingency load-shedding limit constraint based on system operator policy. The proposed methodology is applied to a range of test systems containing up to 10,000 buses with a computational time of up to 3375 s for all 12,706 branch contingencies. Calculating contingency power flows takes 1.3\% of the total solution time using the proposed methodology exploiting the inverse matrix modification lemma.
\end{abstract}


\section{Introduction}
Reliability studies have computational bottlenecks in the fast-growing power systems where the system's state changes faster due to higher penetration of renewable energy, new loads, and more frequent extreme weather. The system is operated closer to the equipment's limits due to increased power needs from the energy transition and electrification, while power system investments are not keeping up with demand growth. Consequently, it is vital to have accurate solutions to ensure the reliability of electricity supply for customers at an affordable price; society is more dependent on it than ever before.

In a traditional reliability management approach, reliability can be accomplished through grid development solutions that consistently meet a deterministic reliability criterion. Operators often use preventive rescheduling of the economic dispatch to meet the criterion. However, such solutions may be prohibitively expensive and not socio-economically efficient when not including corrective actions (generator redispatch, load shedding, and other \ac{TSO} actions after a contingency). More probabilistic and risk-based methods, where reliability and costs are balanced, have been proposed in literature; an example is from the EU-funded project GARPUR: ``\textit{taking into account probabilities of contingencies and using a compound socio-economic objective function blending the costs of \ac{TSO} preventive and corrective actions with a monetization of the risk of service interruptions}''~\cite{garpur_generally_2017}. It is shown with a socio-economic cost assessment that reducing the power flow on critical lines increase power market costs, while interruption costs decrease~\cite{sperstad_probabilistic_2018}. 

\subsection{Related work}
\Ac{SCOPF} is an extensively used tool in reliability studies, whose output could be used to obtain quantitative reliability metrics~\cite{garpur_generally_2017, capitanescu_critical_2016}. The seminal works on \ac{SCOPF} were first published by Alsac and Stott in~\cite{alsac1974optimal} using preventive actions (\acsu{P-SCOPF}) and by Monticelli, Pereira and Granville in~\cite{Monticelli1987Security} using both preventive and corrective actions (\acsu{C-SCOPF}), and since then variants of \ac{SCOPF} has found several applications in the planning and operation of power systems.  
However, accounting for contingencies $\mathcal{C}$ increases the problem size dramatically. Each contingency leads to a different post-contingency state, which ideally is modeled with the same detail as the pre-contingency system. That results in a problem size $|\mathcal{C}|$ times bigger than without accounting for contingencies. Furthermore, if corrective post-contingency actions are available, the number of decision variables also multiplies by $|\mathcal{C}|$. This is a computational challenge, especially with a large system size. In addition, risk and chance constraints could be added to account for uncertainties in demand, generation, or outages, which also increase the complexity of the problem~\cite{Wang2016, ROALD201566}. 

Corrective actions in probabilistic methods need new types of statistical data compared to deterministic methods. These data points are contingency probabilities and the failure rates of the different types of corrective actions. Increasingly good data is available for the former, though it is hard to come by any data for the latter. Corrective action failures could be important to model~\cite{Karangelos2018}, though the present paper will not discuss this. However, the present paper will study the impact of contingency probabilities. Moreover, how the post-contingency system is modeled impacts the outcome. A post-contingency state could be set immediately, minutes, or hours after the contingency, all with different operating limits and possible control actions. At least two post-contingency states are proposed to maintain good security~\cite{Capitanescu2007improving}.

Recent research shows the need for better methods for reliability studies to deal with the computational challenges, as a direct approach of \ac{SCOPF} of large power systems is computationally intractable~\cite{garpur_generally_2017}. 
There is a lot of variation in formulations of \ac{SCOPF} within the literature, with many using an iterative solution procedure. Ensuring post-contingency feasibility while modeling the primary response of generators in a DC \ac{P-SCOPF} is shown for a large power system~\cite{Velloso2021exact}; the primary response is here modeled with binary variables. Corrective actions and outages could be modeled as nodal injections for use in \ac{C-SCOPF}~\cite{Martinez-Lacanina2014corrective}. A probabilistic \ac{C-SCOPF} is proposed in~\cite{he2010optimising}, referred to as an \ac{OPS} approach, which uses an iterative procedure with a particle swarm solution method on a small system. In that approach, weather conditions were incorporated. Analyzing a large power system is possible using a risk-based DC \ac{C-SCOPF} with Lagrangian relaxation and Benders decomposition~\cite{Wang2016}. Several works identify the necessary and sufficient conditions to describe the feasible set of solutions for a DC \ac{P-SCOPF} using pre-processing of the problem, significantly reducing the size of the problem at the expense of the pre-processing time~\cite{Ardakani2013Identification, Weinhold2020Fast}. 
To cope with the computational challenges, machine learning has been employed in reliability management research to predict the system state and generation schedule~\cite{CHEN2022learning}. After a longer initial training of the machine learning model, such methods use a small fraction of the time to find a solution compared to the original \ac{SCOPF} optimization. In the last few years, they have seen great improvement in accuracy and variation in their use within \ac{SCOPF} formulations. Machine learning will not be further discussed in the present paper.
There is a lack of scalable methods of a probabilistic \ac{C-SCOPF} type which would aid operators in system security operational planning optimizing for socio-economic cost benefit. 

\subsection{Contributions and paper organization}
Combining excellent previous work on mathematical methods could help balance computational tractability and accuracy of \ac{SCOPF} and form the motivation for the propositions made in the present paper.
1)~The \ac{IMML}\footnote{Also called the Sherman-Morrison formula, but has been published separately by many authors, so we use the name IMML.} efficiently computes the inverse of a rank-one update of a matrix using the inverse of the original matrix~\cite{sherman1950adjustment}. The present paper uses the \ac{IMML} to solve the contingency power flow problem~\cite{alsac1983sparsity}. 
2)~The sparse matrix solver KLU proposed in~\cite{Davis2010KLU} is used as an alternative calculation method for the contingency power flow. KLU is a linear system solver based on LU-factorization and specifically made to exploit the structure of power system matrices.
3)~The \ac{SCOPF} problem is decomposed into a main problem to which cuts are added, and sub-problems are analyzed using the \ac{IMML} or KLU solver. The seminal work on \ac{SCOPF} used this approach to divide the problem based on the contingency cases~\cite{Monticelli1987Security}. 

The present paper contributes towards reliability studies using a socio-economic \ac{C-SCOPF}, which has stringent requirements for solution speed. Based on the output from the proposed methodology, it is possible to obtain quantitative reliability metrics; the subsequent reliability assessment is not in the scope of the present paper. An efficient and scalable methodology to identify the optimal socio-economic operation of a power system with preventive and corrective actions is proposed. The cutting-plane method is used to iteratively add control actions to the main formulation. The linear cuts are calculated using the \ac{IMML} and the sparse matrix solver KLU. Preventive actions are applied to the pre-contingency (base case) generation schedule, and corrective actions are planned and accounted for, even when the system separates into islands. The proposal extends a previously developed methodology for DC \ac{C-SCOPF}~\cite{Vistnes2023Solving}. The DC approximation is chosen as it is linear and the scope of the study is limited to managing the thermal limits of power system components. The \ac{C-SCOPF} model could be a part of the daily operational security planning of a \ac{TSO}. Before running the \ac{C-SCOPF}, a unit commitment and the latest state estimation of the power system need to be available. In practical applications, further post-processing using AC power flow is recommended to ensure voltage stability and a feasible operation state.

The unique contributions of the present paper are: 
\begin{enumerate}
    \item an improved algorithm based on the previous paper~\cite{Vistnes2023Solving} by the authors, making a)~the method's optimal value of the objective function equal to that of the direct approach, and b)~the solution time faster,
    \item a comparison of the methods for contingency power flow calculations used in the proposed methodology, and 
    \item a detailed case study discussing the two post-contingency states proposed for more accurate steady-state simulations, where combinations of post-contingency states are compared to their \ac{EOC}, \ac{CENS}, and run time.
\end{enumerate}

Next, we discuss the post-contingency modeling decisions in Section \ref{sec:contingency}. Then, Section \ref{sec:method} presents the \ac{C-SCOPF} formulation, contingency power flow solution methods, and the proposed methodology. Further, Section \ref{sec:case_study} presents and discusses the case study. Finally, Section \ref{sec:conclusion} presents a concluding summary of the work.

\section{Sequential post-contingency states and actions}\label{sec:contingency}
\Ac{CENS} serves an active role in probabilistic methods. Its costs and power system impact are weighted against that of generator rescheduling and other possible corrective actions. The flexible load in power systems is increasing, but other loads have grown more reliant on continuous power. There are also soft constraints on what the public tolerates of load shedding before the operator's reputation is impacted. The Norwegian TSO (Statnett) currently has an operating policy that restricts load shedding after a single credible contingency to a maximum of 500 MW (restrictions apply)~\cite{statnett_statnetts_2019}. This is about 2\% of the maximum load in the area of responsibility. A robust constraint is set up based on this to restrict high load-shedding percentages after a contingency. In the case study, the limit is set to 2\%.

Another aspect of post-contingency actions is the temporal evolution. There is an inherent danger in working with static models in that dynamic phenomena can occur between steady-state solutions, leading to unstable trajectories, tripping of relays, or other control mechanisms that cascade the failures. This issue has been tackled in different ways: a line loading risk factor \cite{ROALD201566}, a line failure rate constraint \cite{Subramanyam2022failure}, a robust convex restriction on the load uncertainty set \cite{Dongchan2021robust}, and two post-contingency states~\cite{Capitanescu2007improving}. Focusing on the last, which provides safer control strategies than with one post-contingency state by better modeling the equipment operating limits using the \ac{C-SCOPF}~\cite{Capitanescu2007improving}. 
The model used in the present paper adopts this. Operators can take different actions in the two post-contingency states: the first post-contingency state is set after circuit breakers have tripped, and minutes later, after the reserves are active, the second post-contingency state is set. 

The post-contingency model is a simplification of the balancing reserve market as it is organized in the Nordic countries with a short- and a long-term post-contingency state. The Nordic reserve market is divided into production rescheduling (mitigating problems such as overloaded branches and over- or under-voltage) and frequency reserves (after contingencies such as an outage of production plants or system split). Rescheduling is mainly done through special regulation, where the operator approves a production bid outside the merit order. This can be used both as a preventive action and a corrective action. The \ac{FFR} market is not modeled in the steady-state model used in the present paper. Frequency deviations in the short-term post-contingency state are managed by the \ac{FCR} and the \ac{aFRR} markets. These markets balance the system within its emergency limits, usually valid for 15 minutes of operation. The long-term post-contingency state model results from the \ac{mFRR} market. Compared to markets where the reserves are automatically activated, it has manual activation and keeps the system frequency at its nominal value while allowing long-term emergency ratings (usually valid for 4 hours). These markets, and thus post-contingency states, have different goals; the fast products require fast ramping but minimal energy while the \acl{mFRR} have less stringent ramping but more energy needs to fulfill. Moreover, the markets are cleared separately, and a power plant cannot necessarily provide products for all markets. Hence, the pre-contingency state is a better basis for the market in each post-contingency state than using the previous post-contingency state.

\section{Method}\label{sec:method}
The methodology described in the present paper is a continuation of the work in~\cite{Vistnes2023Solving}. Details of the changes made are provided in this section, while a short summary follows. In the \ac{C-SCOPF} model, a maximum load-shedding constraint is added. Calculation of the post-contingency \ac{PTDF} matrix is done one vector at a time. The iteration procedure is changed and a pre-interation procedure is added.

The power system consists of sets of buses $\mathcal{N}$, branches $\mathcal{B}$, generators $\mathcal{G}$, and demands $\mathcal{D}$, where $N = |\mathcal{N}|$ and $B = |\mathcal{B}|$. Contingencies are denoted $c$ and taken from the contingency set $\mathcal{K}$. Matrix and vector variables are in bold.

\subsection{Corrective Security-Constrained Optimal Power Flow}\label{ssec:SCOPF}
Extending the formulation in~\cite{Capitanescu2007improving}, we use the short-term limits, load shedding and decrease in active power production as control variables, do not only restrict the long-term operational limits to the pre-contingency limits, and add a maximum load shedding constraint. 
\begin{subequations}\label{eq:SCOPF}
\begin{align}
    \min_{\bm{x}^0,\bm{u}^0,\bm{x}^k,\bm{u}^k} \quad &\sum_{g \in \mathcal{G}} c_g p_g + \sum_{k \in \mathcal{K}} \sum_{i \in \bm{u}^{k}} \pi^k u^{k}_{i}
    \label{eq:SCOPF_f}\\
    s.t. \quad & \bm{g^0}(\bm{x}^0,\bm{u}^0) \leq \bm{\overline h}^0  \label{eq:SCOPF_g0}\\
    & \bm{g^k}(\bm{x}^{k},\bm{u}^{k1}) \leq \bm{\overline h}^{k1} &\quad k \in \mathcal{K}   \label{eq:SCOPF_gst}\\
    & |\bm{u}^{k1} - \bm{u}^0| \leq \Delta \bm{u}^{k1} &\quad k \in \mathcal{K} \\
    & \bm{g^k}(\bm{x}^{k},\bm{u}^{k2}) \leq \bm{\overline h}^{k2} &\quad k \in \mathcal{K}  \label{eq:SCOPF_glt}\\
    & |\bm{u}^{k2} - \bm{u}^{0}| \leq \Delta \bm{u}^{k2} &\quad k \in \mathcal{K} \label{eq:SCOPF_delta} 
\end{align}
\end{subequations}
where $\bm{x}$ is the state variables, $\bm{u}$ is the control variables, $\bm{x}^k$, and $\bm{u}^k$, $c_g$ is the marginal cost of generator $g$, $\pi^k$ is the probability of contingency $k$, $\bm{g}$ is abstract notation for the operating limits (including the power flow equations and \ac{TSO} constraints), and $\Delta \bm{u}^k$ is the maximum allowed change of control variables (active power generation and active power decrease or increase between the pre- and post-contingency).
Superscripts $0$ and $k$ represent the base and post-contingency states, respectively. $\bm{\overline h}^0$, $\bm{\overline h}^{k1}$, and $\bm{\overline h}^{k2}$ are operating limits for normal operation, short-term limits after a contingency, and long-term limits after a contingency, respectively.

The modeling of the operating limits $\bm{g}$ in \eqref{eq:SCOPF} is given in the following. The power flow equations are modeled using \acp{PTDF} $\bm{\varphi}$ with the appropriate variables for the specific state. In addition, power balance is needed in all states.
\begin{gather}\label{eq:SCOPF_phi}
    -\bm{\overline h}_l \leq \bm{\varphi} \cdot \left[\bm{u}_g - \left(\bm{\overline h}_d - \bm{u_d}\right)\right] \leq \bm{\overline h}_l \\
    \sum_{g \in \mathcal{G}} u_g - \sum_{d \in \mathcal{D}} u_d = 0
\end{gather}
The post-contingency generation control and load shedding is limited by
\begin{align}\label{eq:SCOPF_pg}
    0 \leq u_{g}^{0} + \Delta u_{g}^{k} \leq \overline h_{g}^{0} &\qquad k \in \mathcal{K}, \quad g \in \mathcal{G} \\
    0 \leq u_{d}^{0} + u_{d}^{k} \leq \overline h_{d}^{0} &\qquad k \in \mathcal{K}, \quad d \in \mathcal{D} \label{eq:SCOPF_pd}
\end{align}
ensuring that active power production is within the generator limits. In the pre-contingency state, $\Delta u_g^k$ and $u_d^k$ are removed from \eqref{eq:SCOPF_pg} and \eqref{eq:SCOPF_pd}.

Further, the maximum allowed post-contingency load shedding $\Gamma$ is modeled as a \ac{TSO} constraint (from the Norwegian \ac{TSO} policy~\cite{statnett_statnetts_2019}). We use this method of load shed constraint as a means to ensure that no outages have higher consequences than are accepted by the operator. It is a form of risk-averse control, a risk policy, for the operator as the linear \ac{VOLL} is not able to capture all costs in extreme circumstances.
\begin{equation}\label{eq:SCOPF_ls}
    \sum_{d \in \mathcal{D}} u_{d} \leq \Gamma
\end{equation}

In summary, the model minimizes the \ac{EOC}, including preventive and corrective actions, such as pre-contingency active power production, active power production re-dispatch, and load-shedding using the probability for each contingency in the objective function~\cite{garpur_generally_2017}. The generator cost is a piecewise linear representation of the generator cost curve, where each piece is modeled as a separate generating unit~\cite{stott1978powerII}. The optimality problem is constrained with all N-1 branch contingencies (transmission lines and transformers). The operating limits are maximum active power in branches (short- and long-term), power balance, generator maximum active power output, generator ramping limit (up and down), and maximum load shedding described by \eqref{eq:SCOPF_g0} -- \eqref{eq:SCOPF_delta} with details in \eqref{eq:SCOPF_phi} -- \eqref{eq:SCOPF_ls}. Suitable operating limits are selected in the pre- and post-contingency states.

A DC power flow neglects reactive power and voltage magnitudes to make the formulation linear. The recent development of convex relaxations of AC power flow could be used to include reactive power and voltage magnitudes. Alternatively, active power losses could be included with a compensation scheme. The present paper does not include power losses, as other parts of the model and solution methodology are highlighted. However, a compensation scheme will be needed in an operational setting. This could easily be included. 

\subsection{Contingency screening}
Modeling the pre- and post-contingency power flow uses the DC power flow approximation. The general modeling of the power flow is first given before a more scalable solution method is shown using the \ac{IMML}.

\subsubsection{DC power flow}
The DC power flow equation for node injected active power $\bm{P}$ is
\begin{gather}\label{eq:dcpf}
    \bm{P} = \bm{\Phi}^\mathsf{T} \cdot \bm{\Psi} \cdot \bm{\Phi} \cdot \bm{\theta} = \bm{H} \cdot \bm{\theta} \\
    \bm{\theta} = \bm{H} \setminus \bm{P} = \bm{X} \cdot \bm{P}
\end{gather}
where $\bm{\Phi}$ is the connectivity matrix, $\bm{\Psi}$ is the diagonal susceptance matrix, $\bm{\theta}$ is the node voltage angles, $\bm{H}$ is the susceptance matrix, and $\bm{X}$ its inverse. The power flow $\bm{F}$ can be found by
\begin{gather}\label{eq:flow_norm_theta}
    \bm{F} = \bm{\Psi} \cdot \bm{\Phi} \cdot \bm{\theta} \\
    \bm{F} = \bm{\Psi} \cdot \bm{\Phi} \cdot \bm{X} \cdot \bm{P} = \bm{\varphi} \cdot \bm{P} \label{eq:flow_norm_ptdf}
\end{gather} 
where $\bm{\varphi}$ is the \ac{PTDF} matrix.

For each contingency, the power flow $\bm{F}^k$ is calculated using the latest iteration of the generation schedule. For a contingency denoted by superscript $k$:
\begin{gather}\label{eq:flow_cont_theta}
    \bm{F}^k = \bm{\Psi}^k \cdot \bm{\Phi}^k \cdot \bm{\theta}^k \\
    \bm{F}^k = \bm{\varphi}^k \cdot (\bm{P} + \Delta \bm{P}^k) =  \bm{\varphi}^k \cdot \bm{P}^k \label{eq:flow_cont_ptdf}
\end{gather}
where $\bm{\varphi}^{k}$ is the contingency \ac{PTDF} matrix (effectively the \ac{LODF}), $\bm{\Psi}^k$ and $\bm{\Phi}^k$ are simple two- and one-element-per-branch contingency changes to $\bm{\Psi}$ and $\bm{\Phi}$, and $\bm{P}^k$ is received from the generation schedule. $\Delta \bm{P}^k$ is the manifestation of the corrective actions. It can often be zero, but is essential to be included when not zero.

Both the \ac{IMML} and KLU solver can be used to calculate $\bm{\theta}^k$ and $\bm{\varphi}^k$. The main differences are highlighted in Table~\ref{tab:methods}. Using a diakoptic solution method, the \ac{IMML} could also be used with system splitting when more than one branch is outaged. However, that solution is not necessarily faster than using a sparse matrix solver~\cite{alsac1983sparsity}.

\begin{table}[tb]
    \caption{DC contingency power flow methods. An ´X' in $\Delta \bm{P}$ (or System Separation) means the method can always be used when the power changes (or the system separates into two or more subsystems) after a contingency. }
    \label{tab:methods}
    \centering
    \begin{tabular}{ccccc} \toprule
        Method & Output & Based on & $\Delta \bm{P}$ & System Separation \\ \midrule 
        I & $\bm{\theta}^k$ & IMML & * & ** \\
        II & $\bm{\theta}^k$ & KLU & X & X \\
        III & $\bm{\varphi}^k$ & IMML & X & ** \\
        IV & $\bm{\varphi}^k$ & KLU & X & X \\ \bottomrule
        \multicolumn{5}{l}{*With pre-processing.} \\
        \multicolumn{5}{l}{**Only implemented for single outage contingencies.}
    \end{tabular}
\end{table}

Method II employs the KLU solver to calculate the contingency power flow using $\bm{H}^{k}$, which is easily found from $\bm{H}$ by reversing the impact of the branch contingency on two diagonal and two off-diagonal elements. 
\begin{gather}\label{eq:M2}
    \bm{\theta}^{k} = \bm{H}^{k} \setminus \bm{P}^k
\end{gather}
The power flow is then obtained using \eqref{eq:flow_cont_theta}. 

Similarly, Method IV employs KLU using $\bm{H}^{k}$ to calculate $\bm{\varphi}^k$ using $\bm{\Psi}^k$ and $\bm{\Phi}^k$.
\begin{gather}\label{eq:M4}
    \bm{\varphi}^{k}[m,:] = \bm{H}^{k} \setminus (\bm{\Phi}^k[m,:] \cdot \bm{\Psi}^k)
\end{gather}
\mbox{$\bm{\varphi}[i,:]$} and $\bm{\varphi}[:,j]$ are notations for the $i$\textsuperscript{th} row and $j$\textsuperscript{th} column of matrix $\bm{\varphi}$, respectively. 
The resulting vector $\bm{\varphi}^k[i,:]$ models the power flow on branch $m$ after the branch outage described by $k$. To calculate $\bm{\varphi}^k$ \eqref{eq:M4} is repeated for all branches or $\bm{\phi}^k[m,:]$ is replaced by $\bm{\phi}^k$ in \eqref{eq:M4}.
Note that $\bm{\varphi}^k$ could later be used to calculate the power flow $\bm{F}^k$ using \eqref{eq:flow_cont_ptdf}, but $\bm{\theta}^k$ cannot be used to find $\bm{\varphi}^k$. 

\subsubsection{Inverse matrix modification lemma}\label{ssec:sm}
The contingency screening process does calculate the post-contingency power flow repeatedly. Thus, efficient performance and scalability are needed for the calculations. Methods I and III adopt the \ac{IMML} to calculate the contingency power flow using the pre-contingency system matrices to achieve the same values as Methods II and IV, respectively. Method I and III (initially) assume a connected system after the contingency. 

Derivation of the equations are found in~\cite{Woodbury1950, Vistnes2023Solving}; Method I calculate $\bm{F}^k$ of branch $l$ from node $i$ to node $j$ using \ac{IMML}, the following are used
\begin{gather}\label{eq:vmc_sm}
    \bm{\theta}^{k} = \bm{\theta}^0 - 
    \frac{\bm{\delta} \cdot \left(\theta^0[i] - \theta^0[j]\right)}{1/H^0[i,j] + \delta[i] - \delta[j]} \\
    \bm{\delta} = \frac{x_l}{H^0[i,j]} \cdot \left(\bm{X}^0[:,i] - \bm{X}^0[:,j]\right) 
\end{gather}
where $x_l$ is the reactance of branch $l$. 
If $\bm{P}^k$ is different than $\bm{P}^0$, an intermediate $\bm{\theta}'$ must be calculated using 
\begin{gather}
    \bm{\theta}' = \bm{X}^0 \cdot \bm{P}^k
\end{gather}
which replaces $\bm{\theta}^0$ in \eqref{eq:vmc_sm}. 

Method III adopts the \ac{IMML} to calculate $\bm{\varphi}^k$ using pre-contingency system matrices and $\bm{\Psi}^k$ and $\bm{\Phi}^k$. 
The intermediate product $\bm{X}^{k}$, the inverse contingency susceptance matrix, could be found directly by inverting the contingency susceptance matrix, $\bm{X}^{k} = (\bm{H} + \bm{\Delta H})^{-1}$. However, using the \ac{IMML} to find $\bm{X}^{k}$ is more efficient using:
\begin{gather}\label{eq:Xc_sm}
    \bm{X}^{k} = \bm{X}^0 - \frac{\bm{X}^0\cdot \bm{\Phi}^0[i,:]\cdot H^0[i,j]\cdot \bm{\delta}}{1+H^0[i,j]\cdot \bm{\delta}\cdot \bm{\Phi}^0[i,:]} \\
    \bm{\varphi}^{k} = \bm{\Psi}^k \cdot \bm{\Phi}^k \cdot \bm{X}^k \label{eq:phi_c}
\end{gather}
The computation can be further streamlined using the knowledge of $\bm{\Phi}^k$ and $\bm{\Psi}^k$. The connectivity matrix $\bm{\Phi}$ has two values for each branch, 1 on node $i$ and -1 on node $j$. $\bm{\Psi}$ only has values on the diagonal. Using that knowledge and simple arithmetics \eqref{eq:Xc_sm} -- \eqref{eq:phi_c} can be written as (superscript 0 is dropped for simplicity):
\begin{gather} \label{eq:phi_sm_vec}
    \bm{\varphi}^k[m,:] = x_m \cdot (\bm{X}[:,n] - \bm{X}[:,o] - \alpha \cdot (\bm{X}[:,i] - \bm{X}[:,j])) \\
    \alpha = \frac{(X[i,o] - X[i,n]) - (X[j,o] - X[j,n])}{1 / x_l + (X[i,j] - X[i,i]) - (X[j,j] - X[j,i])}
\end{gather}
for an overload on branch $m$ from node $n$ to node $o$. \mbox{Eq.~\eqref{eq:phi_sm_vec}} is more efficient because $\alpha$ is a single number and $\bm{\varphi}^k[m,:]$ is calculated using only a few vectors.

Methods I and III are slightly changed when a single contingency cuts off nodes radially connected to the rest of the system post-contingency. The cut-off nodes are disconnected from the ground in $\bm{H}^0$ by zeroing out the corresponding rows and columns before applying the described methods. For higher-order contingencies, a diakoptic solution method could be adopted~\cite{alsac1983sparsity}.

\subsection{Algorithm}\label{ssec:algorithm}
\begin{figure}[tb]
    \centering
    \includegraphics[width=\linewidth]{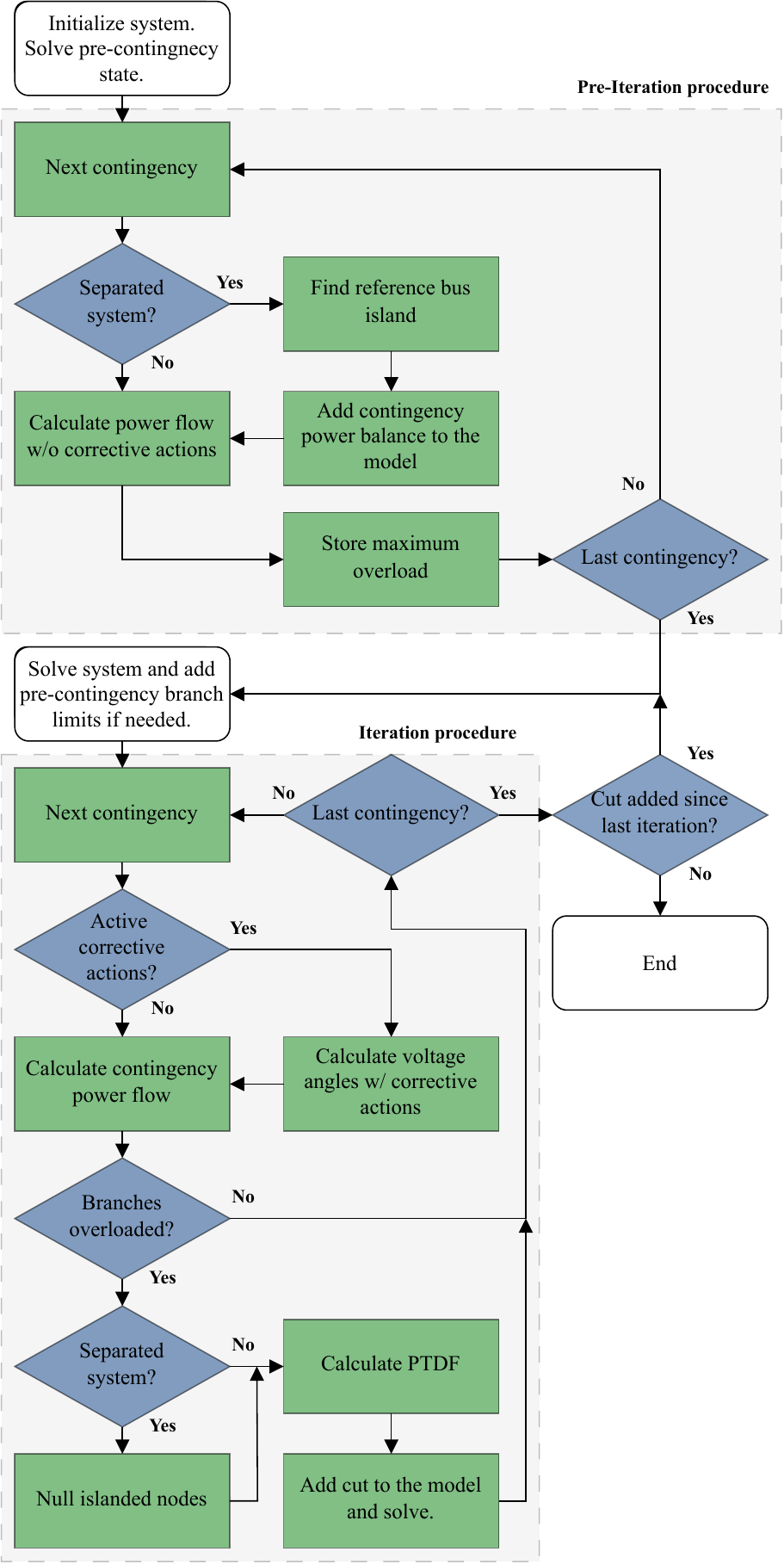}
    \caption{Flow chart of the proposed algorithm for solving a \ac{C-SCOPF}.}
    \label{fig:SCOPF}
\end{figure}

The solution algorithm, an improvement over the previous version~\cite{Vistnes2023Solving}, is shown in Fig.~\ref{fig:SCOPF}. The optimization problem \eqref{eq:SCOPF} is decomposed into a main problem, consisting of the pre-contingency state, and sub-problems, each consisting of a post-contingency state. The total number of system states is $2C+1$. A constraint (cut) is added to the main problem for each post-contingency overloaded branch. Please refer to~\cite{Vistnes2023Solving} for the deduction of equations for the cuts.

A particularly important algorithm step is in handling the power flow constraints. This is changed from \cite{Vistnes2023Solving} and exploits the fast contingency power flow computation. At the program's start, no power flow constraints are added, and the model is similar to an economic dispatch problem, only considering generation dispatch, and disregarding the system connectivity. After the initial generation dispatch is found, branch constraints are added to the pre-contingency state for overloaded branches. Consequently, the optimization problem is much smaller and thus easier and faster to solve than using all branch constraints. 

A pre-iteration procedure is added to the algorithm that serves a twofold purpose. It adds corrective actions for every contingency that results in system separation. The island without the reference node shed all load. Moreover, the pre-iteration procedure does a contingency screening by calculating the contingency power flow using Method I (or Method II for comparison). It ranks the contingencies from maximum branch overload to the lowest (or no overload). During this pre-iteration procedure the \ac{SCOPF} is not re-solved; the system is only re-solved after all variables for corrective actions have been added to the model, as seen in the middle of Fig.~\ref{fig:SCOPF}. As a double benefit, the ranking of contingencies results in cuts being added for the most `difficult' contingencies first; the large changes to the generation schedule were forced earlier in the main iteration procedure, resolving other contingencies with a similar but lower impact on branch loading. Both sorting based on short- and long-term branch limits were tried. Sorting based on short-term limits was deemed better as these limits more often change the pre-contingency schedule to the end optimal schedule in our testing. The model is re-solved right before the main iteration procedure to find the new pre-contingency generation schedule. In addition, the re-solving balances the remaining sub-system using post-contingency corrective rescheduling for each contingency where the system separates.

In the second half of the algorithm (the bottom of Fig.~\ref{fig:SCOPF}) the necessary cuts are added to the model to enforce the branch flow limits $g$ in \eqref{eq:SCOPF}. First, \eqref{eq:SCOPF_g0} is re-evaluated. Then, \eqref{eq:SCOPF_gst} and \eqref{eq:SCOPF_glt} are ensured through the iteration procedure. Intermediate voltage angles are calculated if there are corrective actions for the post-contingency state, before the contingency power flow is calculated. Next, suppose there are overloaded branches in the system after a contingency with the current generation schedule. In that case, $\bm{\varphi}^k$ is calculated by Method III and is used in the addition of a cut for the short- or long-term branch limits. 

As in the pre-iteration procedure, the contingency power flow is always calculated by Method I, a change from the algorithm presented in~\cite{Vistnes2023Solving} resulting from the implementation of system separation of single outages in Method I. Moreover, the full dense matrices ($X$ and $\varphi$) are not needed in the algorithm; only the parts of these matrices corresponding to the overloaded branches are needed for each contingency. Lastly, the standard linear solver in Julia was switched to the KLU solver to improve speed. 

\section{Case Study}\label{sec:case_study}
The efficiency of the proposed methodology is demonstrated by reporting results from a range of systems, from 24 to 10,000 buses~\cite{Subcommittee1979rts, Birchfield2017Synthetic}, with accompanying discussions. The systems are set to peak load according to the system description. \ac{VOLL} for load shedding are found in ~\cite{Firuzabad2000Impact}. The DC approximation models the systems using all buses and branches (lines and transformers). The probability of a contingency $\pi^k$ is set as the average yearly outage probability for an hour for each branch. All N-1 branch contingencies are considered in each \ac{SCOPF} for each system. For the 10,000 bus system, considering all 12,706 branch contingencies using the formulated \ac{C-SCOPF}, there are 85 million decision variables and 373 million constraints. In other terms, this formulation is intractable for direct solution. However, using the proposed algorithm, the problem starts out with 6655 decision variables and 16,655 constraints, iteratively growing as needed but never approaching the size of the naive problem formulation.

Simulations are run for three variants of \ac{SCOPF}: 
\begin{itemize}
  \item SCOPF: Only a pre-contingency state is optimized.
  \item P-SCOPF: A pre-contingency and one post-contingency state are optimized. Post-contingency overloads are based on the long-term limits by using preventive control actions.
  \item C-SCOPF: A pre-contingency and two post-contingency states (both short and long-term limits) with preventive and corrective control actions are optimized.
\end{itemize}

All results are gathered on an Intel\textsuperscript{\textregistered} Core\textsuperscript{TM} i9-9920X, 128 GB RAM server. The methodology is implemented in the Julia programming language v1.9~\cite{Bezanson2017Julia} using JuMP v1.22.2~\cite{Lubin2023JuMP} and PowerSystems v3.3.0~\cite{Lara2021PowerSystems}, linear solver KLU v0.6.0~\cite{Davis2010KLU}, and problem optimization using the Gurobi v10.0 solver~\cite{gurobi}. The reported solution times include the model preparation time, intermediate calculations, and LP solver time. 

The case study is divided into two parts. First, the computational time of the proposed solution algorithm is shown for a breadth of different system sizes. Further, feasibility and the improved computational time for large power systems enable analysis of the systems using \ac{C-SCOPF}. A sensitivity analysis is presented in the second part of the case study.

\subsection{Computational time}
\subsubsection{DC power flow}
\begin{figure}[tb]
    \centering
    \includegraphics[trim={2.7cm 5.cm 2.3cm 4.5cm}, clip, width=\linewidth]{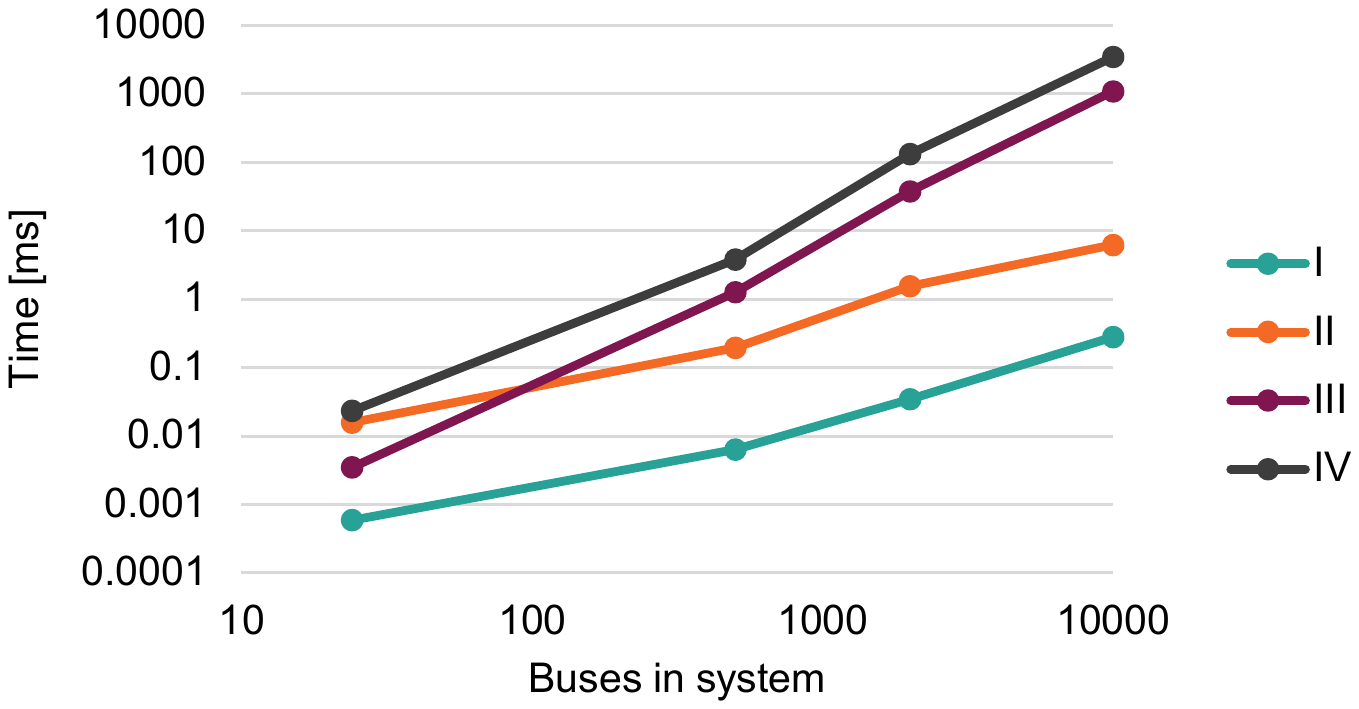}
    \caption{Computational time of DC contingency power flow with change of post-contingency power $\Delta P$ and system separation.}
    \label{fig:time_dcpf}
\end{figure}

Fig.~\ref{fig:time_dcpf} shows the computational time for DC contingency power flow. Method IV is the slowest method for all system sizes. The difference between methods increases with the number of nodes in the system. For the 10,000 bus system, Method III is three times faster than Method IV. Further, Method II is two orders of magnitude faster, and Method I is over a magnitude faster again. In total, Method I using \ac{IMML} is four orders of magnitude faster than Method IV based on KLU and $\bm{\varphi}^{k}$.

$\bm{X}$ and $\bm{\varphi}$ are typically dense matrices for all power systems. They are derived from the $\bm{B}$, $\bm{\Phi}$, and $\bm{\Psi}$ matrices that typically for large power systems are very sparse. Moreover, the computational power and memory needed to calculate and store the dense matrices can be considerable. Therefore, if only part of a matrix is needed, that is the only part that should be calculated. Only one row at a time from $\bm{\varphi}^{k}$ is needed when iteratively adding cuts. Using Method III and Method IV for calculation of the $\bm{\varphi}^{k}$ sub-matrix (a vector) leads to similar computational times as for Method I and Method II calculating the power flow, respectively.

\subsubsection{SCOPF}
The proposed algorithm is used on system sizes up to 10,000 buses to show the algorithm's scalability.  
Table~\ref{tab:runtime} shows the computational time using different system sizes and \ac{SCOPF} variants. The computational time (time to solution) for the 10,000 bus system with no contingencies is \SI{0.42}{s}, while considering preventive and corrective security actions for contingencies gives a computational time of \SI{3375}{s}. 

\begin{table}[tb]
    \centering
    \caption{Computational time for solving the specified \ac{SCOPF}.}
    \label{tab:runtime}
    \begin{tabular}{@{}lrrrrr@{}} \toprule
          & & & \multicolumn{3}{c}{Computational time [s]} \\ \cmidrule{4-6}
         System & Nodes & Branches & SCOPF & P-SCOPF & C-SCOPF \\ \midrule
         IEEE RTS & 24 & 38 & $<$0.01 & $<$0.01 & $<$0.01 \\
         ACTIVSg500 & 500 & 597 & 0.01 & 0.11 & 1.33 \\
         ACTIVSg2000 & 2000 & 3206 & 0.05 & 1.73 & 18.85 \\
         ACTIVSg10k & 10,000 & 12,706 & 0.42 & 168 & 3375 \\ \bottomrule
    \end{tabular}
\end{table}

\begin{figure}[tb]
    \centering
    \includegraphics[trim={1.8cm 3.5cm 1.8cm 4.6cm}, clip, width=\linewidth]{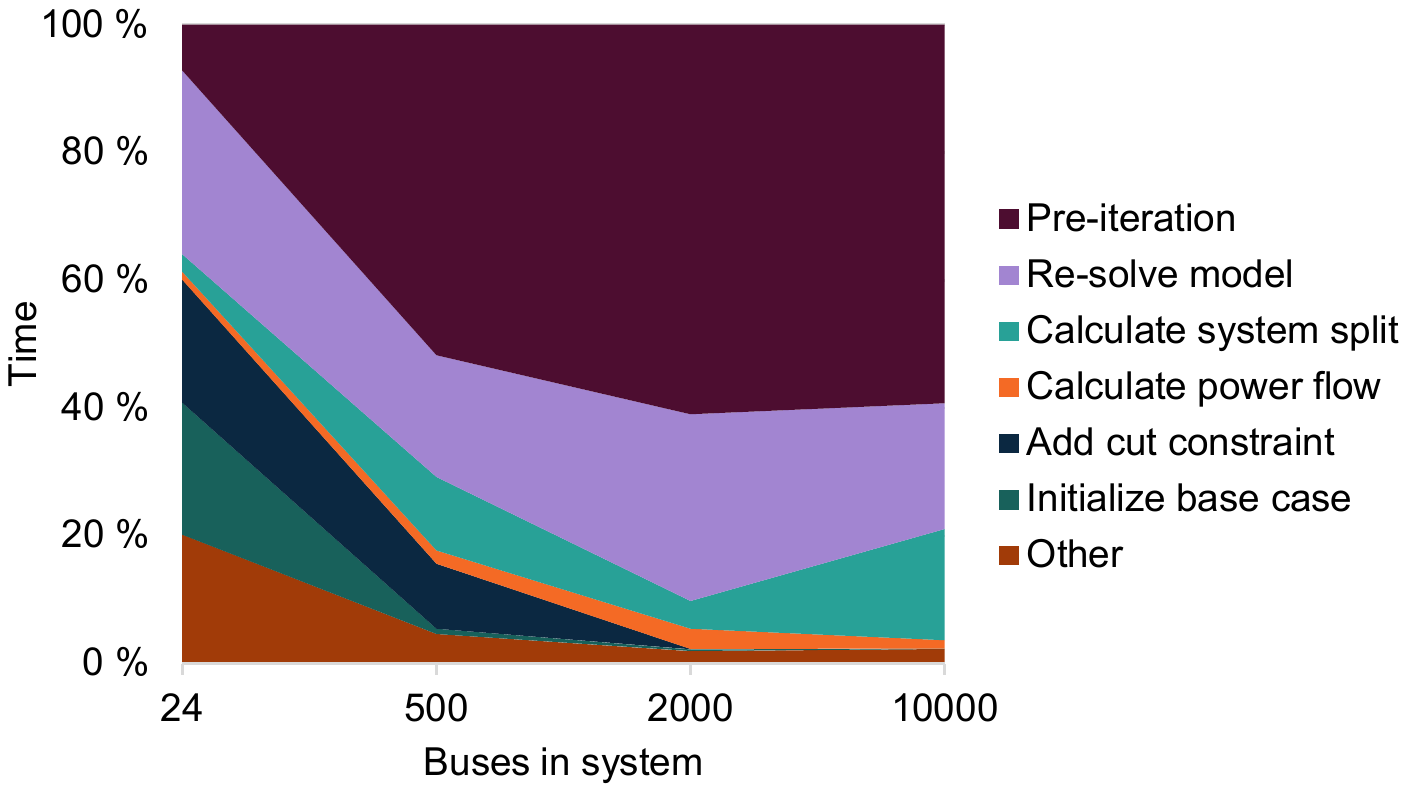}
    \caption{Breakdown of the time used for each part of the solution algorithm.}
    \label{fig:runtime_breakdown}
\end{figure}

A breakdown of the computational time into different categories is shown in Fig~\ref{fig:runtime_breakdown}. As also seen in Table~\ref{tab:runtime}, setting up and solving the problem using the \ac{SCOPF} formulation (initialize pre-contingency in Fig.~\ref{fig:runtime_breakdown}) takes negligible time for the larger systems. This shows that a focus on post-contingency modeling is warranted. Moreover, using better approximations for the pre-contingency could be considered while keeping the DC post-contingency approximation. 
Calculation of contingency power flow, a focus area of the present paper, only uses 0.9--3.2\% of the total computational time. Hence, the contingency screening is not a limiting factor in the algorithm. 19.1--29.1\% of the time is used inside the optimization solver after new cuts are added to the model. 

The pre-contingency procedure takes up to 60.6\% of the computational time for the larger systems. This time is mostly spent adding variables and constraints for the two post-contingency states in the model for the system splitting contingencies. Note that removing the pre-iteration procedure would not eliminate the computation time; it just spreads out during the iteration procedure. Moreover, implementation, choice of optimization solver, solver interface, and programming language influence this block of time greatly. Without these variables one would need another scheme to ensure enough ramping capacity in the correct locations to balance the power. An alternative is to use a predefined distributed slack that ensures reserve capacity on the participating generators for all foreseeable contingencies. This would also greatly reduce the model's size such that the time spent on re-solving the model would also be reduced.

\subsection{Sensitivity analysis}

The sensitivity of \ac{EOC} from the method \ac{C-SCOPF} is analyzed with respect to the input parameters using the synthetic 500-bus SouthCarolina500 system~\cite{Birchfield2017Synthetic}. The system has 597 branch contingencies, of which 254 branch contingencies split the system by isolating one bus.

\begin{figure}[tb]
    \centering
    \includegraphics[trim={0cm 2cm 0.5cm .5cm}, clip, width=\linewidth]{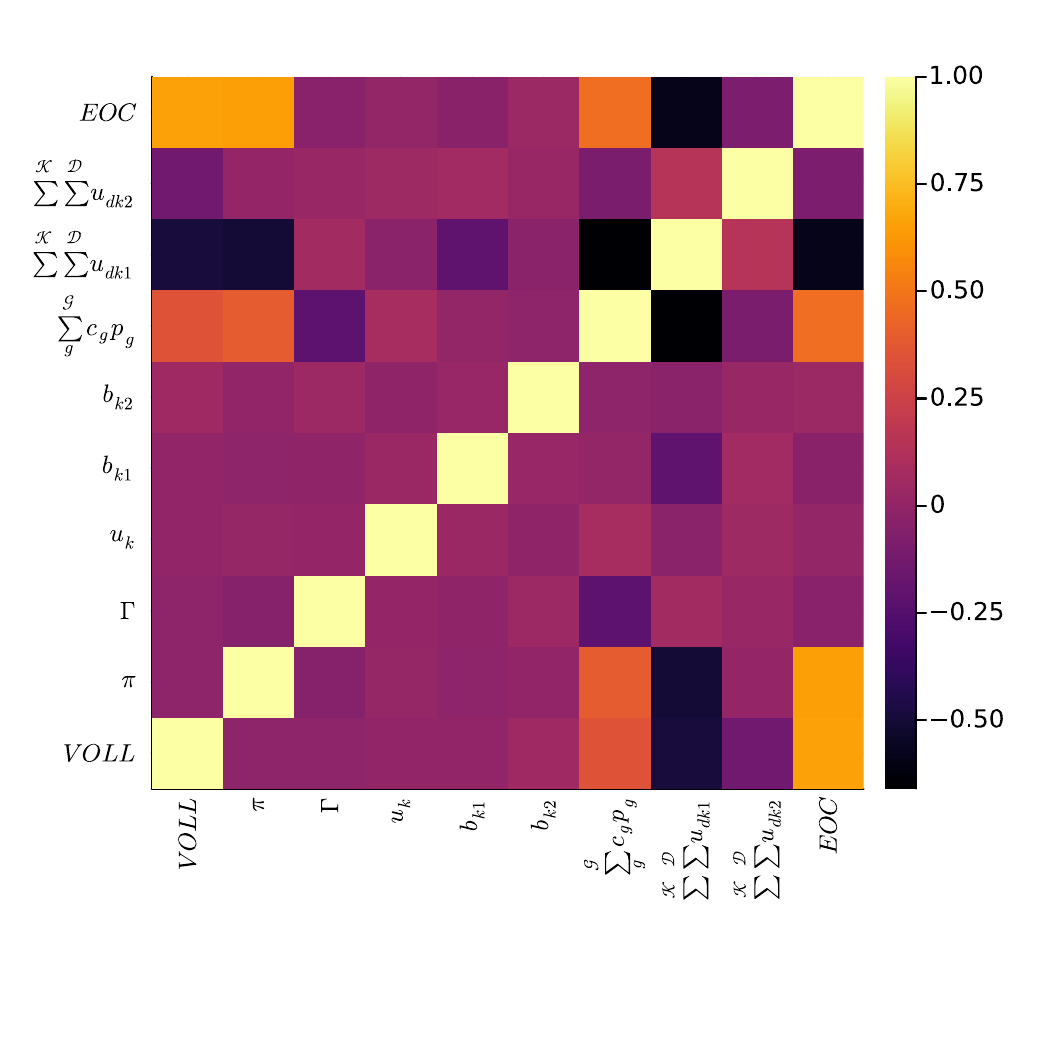}
    \caption{Heatmap for the 500-bus system.}
    \label{fig:500bus_heatmap}
\end{figure}

\begin{figure*}[tb]
    \centering
    \includegraphics[width=\linewidth]{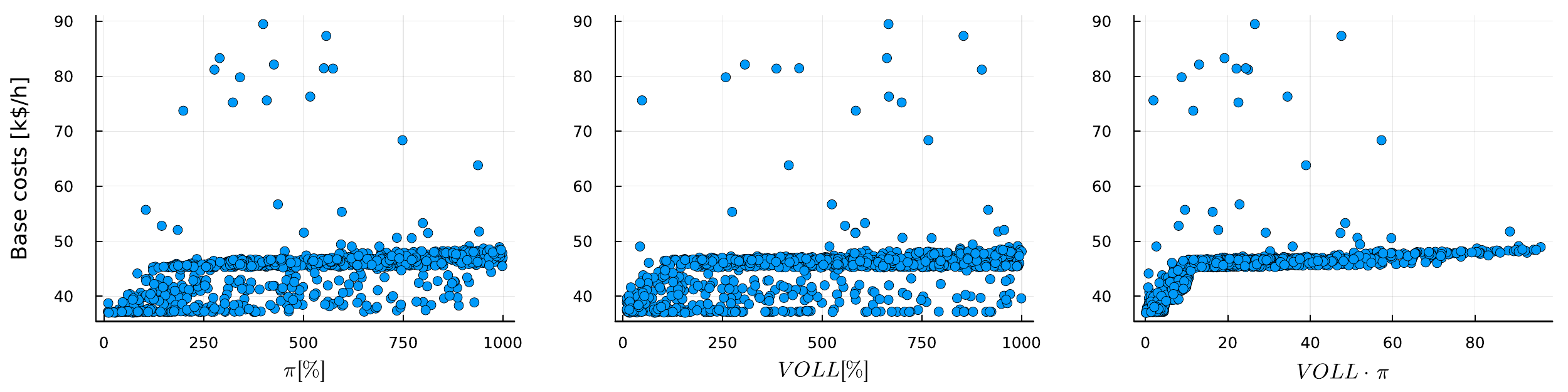}
    \caption{Scatter plot of base costs vs \ac{VOLL} multiplied with $\pi$ for the 500-bus system.}
    \label{fig:pg0_vs_VOLL_p}
\end{figure*}

The sensitivity of the input parameters to the output values is found by random sampling. 1000 random samples were drawn for each parameter (see Table~\ref{tab:parameters}) as input to a \ac{C-SCOPF} of the 500-bus system. The result is summarized in Fig.~\ref{fig:500bus_heatmap} which shows the correlation between the input parameters, \ac{VOLL}, $\pi$, $\Gamma$, $u^k$, $b^{k1}$, and $b^{k2}$, and the output costs, base costs (pre-contingency costs), short-term corrective load shedding, long-term corrective load shedding, and \ac{EOC}. Note that \ac{EOC} is a combination of the other output costs, and it is \ac{ENS} and not \ac{CENS} in the heatmap. The input parameters were drawn from uniform distributions and did not show any correlation in Fig~\ref{fig:500bus_heatmap}. The base costs and short-term \ac{ENS} show a strong negative correlation. While short-term \ac{ENS} also negatively correlates with the other parameters, it is not as strong. In addition, the base costs and the \ac{EOC} positively correlate with $\pi$ and \ac{VOLL}, with \ac{EOC} showing the strongest correlation. Further investigating this relationship, Fig.~\ref{fig:pg0_vs_VOLL_p} shows the two parameters $\pi$ and \ac{VOLL} multiplied. It is evident from looking at Fig.~\ref{fig:pg0_vs_VOLL_p} that both $\pi$ and \ac{VOLL} need to have low values for the \ac{C-SCOPF} to schedule the cheapest pre-contingency generation.
\begin{table}[tb]
    \centering
    \caption{Parameter multipliers with its uniform distribution.}
    \label{tab:parameters}
    \begin{tabular}{cc} \toprule
        Parameter & Range \\ \midrule
        VOLL & 0.1--10 \\
        $\pi$ & 0.1--10 \\
        $\Gamma$ & 0.1--50 \\
        $u^k$ & 1--10 \\
        $\Delta u^{k2}$ & 10 \\
        $b^{k1}$ & 1.25--1.5 \\
        $b^{k2}$ & 1.0--1.25 \\ \bottomrule
    \end{tabular}
\end{table}

\begin{figure}[tb]
    \centering
    \includegraphics[width=\linewidth]{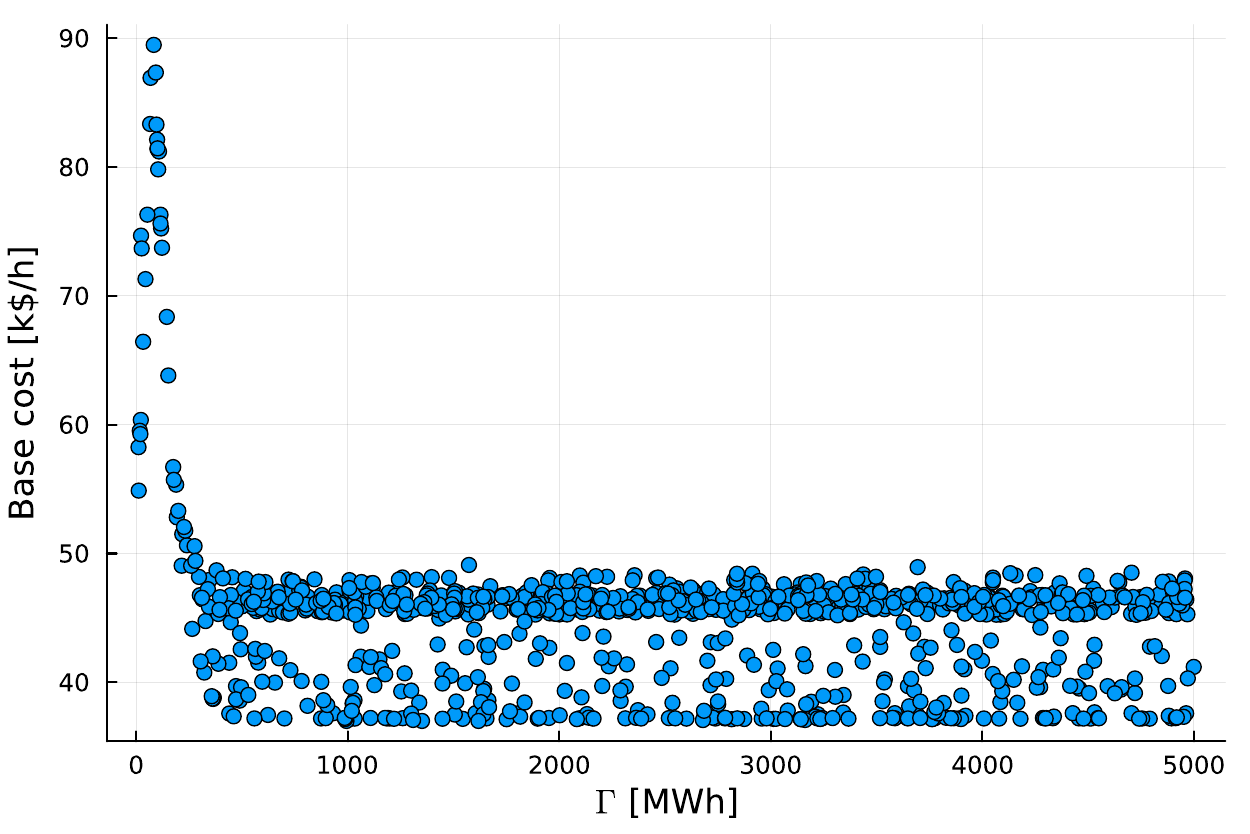}
    \caption{Scatter plot of base costs vs $\Gamma$ for the 500-bus system.}
    \label{fig:500bus_gamma}
\end{figure}

Maximum allowed load shed per contingency $\Gamma$ has a slight negative correlation with base costs, as shown in Fig~\ref{fig:500bus_heatmap}. The scatter plot in Fig.~\ref{fig:500bus_gamma} shows the $\Gamma$ and base cost correlation is limited; for low values of $\Gamma$ there is a strong negative correlation. However, for higher values the base cost is not affected by $\Gamma$.

The sensitivity analysis suggests that the main drivers of the input parameters for low base cost are low $\pi$ and \ac{VOLL}, while high base cost is driven by the restriction of post-contingency load shedding ($\Gamma$). Increasing the generator post-contingency rescheduling cost has less impact on the base costs in this case study of the 500-bus system. 

\begin{figure}[tb]
    \centering
    \includegraphics[width=\linewidth]{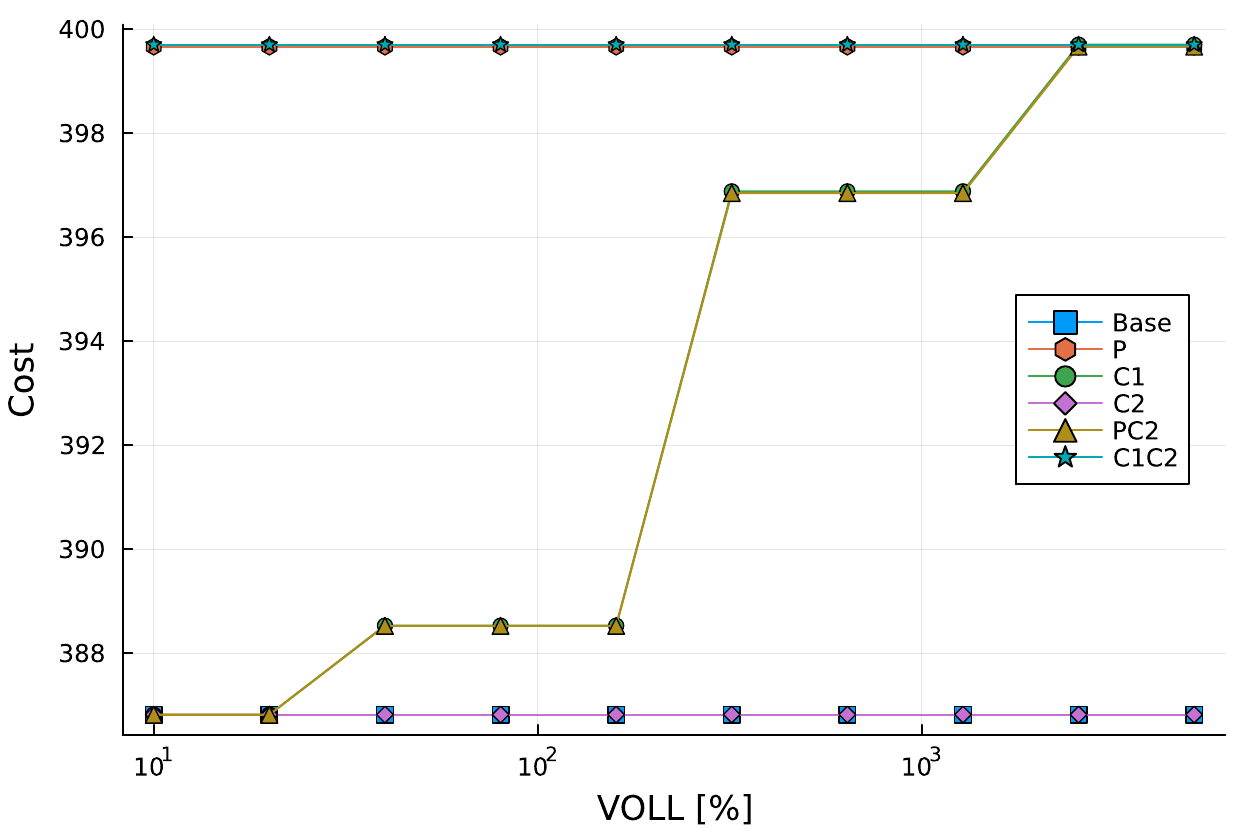}
    \caption{Base case cost for 500-bus system changing the \ac{VOLL}. Acronyms are for the pre-contingency state (P), and for the post-contingency states they are the short- (C1) and long-term (C2). }
    \label{fig:cens}
\end{figure}

In the \ac{C-SCOPF} objective the generation cost is balanced against \ac{CENS}, both pre- and post-contingency, where \ac{CENS} is \ac{VOLL} times \ac{ENS}. Thus, \ac{VOLL} is as essential as the marginal cost of generators in a \ac{SCOPF}. However, the \ac{VOLL} is much harder to quantify than the marginal generation cost as the \ac{ENS} could be from many different customers whose perceptions of \ac{VOLL} differ. Removing all system splitting contingencies, Fig.~\ref{fig:cens} shows how changing the \ac{VOLL} changes the pre-contingency generation cost for the \ac{SCOPF} variants. For up to 30\% \ac{VOLL} the \ac{C-SCOPF} makes the same schedule as the \ac{SCOPF} (without any contingencies). Further, by increasing the \ac{VOLL}, the base cost increases as more preventive measures are scheduled to lessen the \ac{CENS} post-contingency. Above 3000\% \ac{VOLL}, the base cost of the \ac{C-SCOPF} is the same as that from the \ac{P-SCOPF}.

With respect to the discussion in this case study, it is essential to note that these results are gathered on specific test-power systems, and the conclusions are not necessarily general. However, the trends seem similar for the systems in this case study. 

\section{Conclusion}\label{sec:conclusion} 
A fast and scalable iterative solution methodology of a socio-economic \ac{C-SCOPF} with two post-contingency states is proposed. Moreover, essential to good scalability in the methodology are the contingency power flow and optimization cuts, which are calculated using the \ac{IMML}. Results from the case study using the proposed methodology did show the feasibility of a 10,000-bus system \ac{C-SCOPF} on modest computer hardware, an Intel\textsuperscript{\textregistered} Core\textsuperscript{TM} i9-9920X with 128 GB RAM. Still, the method depends on the number of contingencies constraining the system in the optimum, including contingencies that separate the system. Optimizing using only preventive actions for the 10,000-bus system took \SI{168}{s} and also including two-post-contingency states with corrective actions resulted in a computational time of \SI{3375}{s}. There is still room for further improvements, especially regarding the pre-iteration procedure where the post-contingency power balance is ensured. Alternatively, a modeling option is to use distributed slack with pre-contingency reserved capacity for ramping. The sensitivity analysis shows cost-correlations with $\pi$, \ac{VOLL}, and $\Gamma$. $\pi$ are system state estimations and could be erroneous. \ac{VOLL} and $\Gamma$ are set parameter quantities based on customer values and system operator policy, respectively, but could be hard to quantify. These parameters should be further investigated for other systems using \ac{C-SCOPF}.

A natural extension of the proposed method includes multi-contingency outages and islanded operation of the system.
Using the results from the proposed \ac{C-SCOPF} and a \ac{C-SCOPF} with only one post-contingency state, a dynamic time-domain simulation could be conducted to show how the dynamic security changes with one post-contingency state, compared to two post-contingency states.
Integrating methods for contingency screening and concurrent contingency calculation will further increase the applications for which the proposed methodology can be used.

\let\bibfont\footnotesize
\printbibliography

\end{document}